\documentclass[12pt]{article}
\usepackage{latexsym}
\usepackage{amsfonts}
\usepackage{epsfig}
\newcommand{\binom}[2]{\left(\begin{array}{c}#1\\#2\end{array}\right)}
\def\ord{\mathrm{ord\,}}
\def\C{{\mathbf{C}}}
\def\q{\mathbf{q}}
\def\e{\mathbf{e}}
\def\k{\mathbf{k}}
\def\f{\mathbf{f}}
\def\z{\mathbf{z}}
\def\a{\mathbf{a}}
\def\x{\mathbf{x}}
\def\y{\mathbf{y}}
\def\R{{\mathbf{R}}}
\def\bR{{\mathbf{\overline{R}}}}
\def\P{{\mathbf{P}}}
\def\Poly{{\mathrm{Poly}}}
\def\Rat{{\mathrm{Rat}}}
\begin{document}
\title{Elementary proof of the B. and M. Shapiro conjecture for rational functions}
\author{Alex Eremenko
and Andrei Gabrielov\thanks{Both authors were supported by NSF.}}
\date{October 10, 2005}
\maketitle
\begin{abstract}
We give a new elementary proof of the following theorem:
if all critical points of a rational function $g$ belong to the real line
then there exists a fractional linear transformation $\phi$ such that
$\phi\circ g$ is a real rational function.
Then we interpret the result in terms of Fuchsian differential equations
whose general solution is a polynomial and
in terms of electrostatics.
\end{abstract}

One of the many equivalent formulations of the Shapiro conjecture is
the following. Let $\f=(f_1,\ldots,f_p)$ be a vector of polynomials in
one complex variable, and assume that the Wronski determinant
$W(\f)=W(f_1,\ldots,f_p)$ has only real roots. Then there exists a matrix
$A\in {\mathrm{GL}}(p,\C)$ such that $\f A$ is a vector of real polynomials.

This conjecture plays an important role in real enumerative geometry
\cite{Soterag,SotShap}, theory of real algebraic curves \cite{KS}
and has applications to control theory
\cite{RosSot,SIAM}.
There is a substantial numerical evidence \cite{SotShap}
in favor of the conjecture.

In \cite{Ann} we proved the Shapiro conjecture in the first non-trivial
case $p=2$. The proof was quite complicated, and its main drawback
from the point of view of generalizations to higher dimensions was the
use of the Uniformization theorem.

In this paper we give a new proof, not using the Uniformization theorem.
The new proof is also much simpler
than the arguments in \cite{Ann}. 

Consider a non-constant
rational function $g=f_1/f_2$, and assume that the polynomials
$f_1$ and $f_2$ are co-prime. Then the degree of $g$ is given by
$d=\max\{\deg f_1,\deg f_2\}$, and the roots of the Wronski determinant
$W(\f)=f_1f_2^\prime-f_1^\prime f_2$ coincide with the critical
points of $g$.  Let us call two rational functions $g_1$ and $g_2$ 
equivalent if $g_1=\phi\circ g_2$ for some fractional-linear transformation
$\phi$. Evidently, equivalent rational functions have the same critical
points. So our result is
\vspace{.1in}

\noindent
{\bf Theorem 1.} {\em If all critical points
of a rational function are real
then it is equivalent to a real rational function.}
\vspace{.1in}

It is enough to prove this theorem for rational functions with simple
critical points. The general case then follows by a limiting process.

It is known \cite{Goldberg} that for given
$2d-2$ points in the complex plane
in general position, there exist
\begin{equation}
\label{Catalan}
u_d=\frac{1}{d}\binom{2d-2}{d-1},
\end{equation}
the $d$-th Catalan number of classes of rational functions of degree $d$
with these critical points.
It turns out that the general position assumption
in this result can be removed if the given points are real.
Moreover, the following result turns out to be equivalent to Theorem 1:
\vspace{.1in}

\noindent
{\bf Theorem 2.} {\em For any given $2d-2$ distinct
points on the real line, there exist
exactly $u_d$ distinct classes of rational functions of degree
$d$ with these critical points.}
\vspace{.1in}

It follows from Theorem 1 that each of these $u_d$ classes contains
a real function.
The assumption that the critical points are real is essential in
Theorem 2: for $2d-2$ complex points, the number of rational functions
of degree $d$ with these critical points can be less than $u_d$.

Equivalence of theorems 1 and 2 was known for some time, see, for example,
\cite{SotShap}.
 
To state a generalization of Theorem 2 to the case of multiple critical points,
we recall the definition of {\em Kostka numbers.}
Let $\a=(a_1,\ldots,a_q)$ be a vector of integers satisfying
\begin{equation}
\label{inte} 
1\leq a_j\leq d-1,\quad\sum_{j=1}^q a_j=2d-2.
\end{equation}
Consider the Young diagrams of shape $2\times(d-1)$. They consist of two
rows of length $d-1$. A {\em semi-standard Young tableau}
SSYT of shape $2\times(d-1)$ is a filling of such a diagram by positive
integers,
such that an integer $k$ appears $a_k$ times, the entries are strictly
increasing in columns and non-decreasing in rows.
The corresponding Kostka number $K_\a$ is the number of such SSYT.
The number $K_\a$ does not change if the coordinates of $\a$ are permuted
\cite[Thm. 7.10.2]{Stanley}.
\vspace{.1in}

\noindent
{\bf Theorem 3.} {\em For given $\a$ satisfying $(\ref{inte})$,
and given real points $x_1<x_2<\ldots<x_q$, 
there are exactly $K_\a$ classes of rational functions of degree $d$
with critical points at $x_j$ of multiplicity $a_j$.}
\vspace{.1in}

We obtain Theorem 2 as a special case when $q=2d-2$ and all $a_j=1$.
Theorem 3 is true for {\em generic} complex $x_j$; this was derived
by Scherbak \cite{Sch} from the results in \cite{Varsch}.
Theorem 3 was first proved in \cite{egsv},
where a result from \cite{Ann} was used.
We include a proof here to show that it can be achieved with the
same elementary tools as theorems 1 and 2, and no heavy machinery
from \cite{Ann} is needed. 

A self-contained proof of Theorems 1 and 2 is given in Section 1.
In Section 2, we discuss two equivalent reformulations of
these theorems, which in our opinion are of independent interest,
and also might be helpful for proving
the Shapiro Conjecture in higher dimensions. Many other reformulations
of the Shapiro conjecture are contained in \cite{SeS,Soterag,SotShap}. 

We thank Boris Shapiro, Frank Sottile and Alexander Varchenko for their
useful comments on this paper.

\section{Proof of the main theorems}

\subsection{Wronski map}

We recall the necessary facts and definitions.
Let $\f=(f_1,f_2)$
be a pair of linearly independent polynomials of degree at most $d$.
They span a $2$-dimensional subspace in the space 
of all polynomials
of degree at most $d$, and thus define a point in the Grassmannian
$G=G(2,d+1)$. Two pairs of polynomials are equivalent if they span the same
subspace.

The Wronski determinants of equivalent pairs are proportional
non-zero polynomials
of degree at most $2d-2$. Classes of proportionality of such polynomials
form a space $\Poly^{2d-2}$
which can be identified
with the projective space $\P^{2d-2}$. So taking Wronski
determinant defines a map
\begin{equation}
\label{W}
W:G(2,d+1)\to \Poly^{2d-2}
\end{equation}
which is called the {\em Wronski map}.
The real Grassmannian $G_\R$ or the real projective space $\Poly_\R$
consist of those points whose coordinates (coefficients
of the polynomials) can be chosen real.
It is clear that $W$ sends $G_\R$
to $\Poly^{2d-2}_\R$. 
The Wronski map is a finite regular map of compact algebraic manifolds,
and its degree can be defined as the number of preimages of a generic
point. This number turns out to be the Catalan number
$u_d$, see, for example,
\cite{Goldberg}. If the Grassmannian $G$ is embedded in a projective space
by the Pl\"ucker
embedding, then the Wronski map
$W$ becomes a restriction of a linear projection on $G$.
Thus the degree of $W$
is the same as the degree of the Grassmann
variety, that is the number of intersections of a generic
subspace of codimension $2d-2$ with the Pl\"ucker embedding of $G$.  

Using this notation, theorems 1 and 2 can be restated as follows:
\vspace{.1in}

\noindent
1. {\em  The full preimage of
a polynomial with all real roots under the Wronski
map consists of real points in $G$.}
\vspace{.1in}

\noindent
2. {\em Every polynomial in $\Poly^{2d-2}$ with distinct real
roots has exactly $u_d$
distinct preimages under the Wronski map.}
\vspace{.1in}

If a pair $(f_1,f_2)$ represents a point of $G$, then $g=f_1/f_2$
is a non-constant
rational function
of degree at most $d$. If two pairs of polynomials represent the same point
of $G$, then the corresponding rational functions are equivalent.
This defines a map $r$ from $G$
into the set $\Rat^d$ of equivalence classes of non-constant rational
functions of degree at most $d$. This map is not
injective because polynomials in
a pair can have a common factor. More precisely,   
let $Z_0\subset G$ be the locus of points corresponding to
pairs of polynomials having a non-constant
common factor, and $Z_1$ the locus corresponding to polynomials
of degree less than $d$. Put $Z=Z_0\cup Z_1$. Then 
\begin{equation}
\label{r}
r:G\backslash Z_0\to\Rat^d
\end{equation} is a bijection.
The standard topology on $G$ can
be defined as induced by the Pl\"ucker embedding,
and the topology on $\Rat^d$ is of uniform convergence with respect to
the spherical metric.
The map $r$ is continuous on $G\backslash Z$ but
not continuous on the whole $G$.
The following weaker continuity property of $r$ holds on the whole
Grassmannian.
\vspace{.1in}

\noindent
{\bf Proposition 1.} {\em Let $(p_j)$ be a converging sequence in $G$, 
$p_j\in G\backslash Z$ and $p=\lim p_j$ is represented
by a pair of polynomials with a common factor $q$.
Let $z_1,\ldots,z_k$ be the roots of $q$. Then $r(p_j)\to r(p)$
uniformly on compact subsets of $\C\backslash \{ z_1,\ldots,z_k,\infty\}$.
One does not have to include $\infty$ if the degree of $p$ is
the same as that of $p_j$.}
\vspace{.1in}

The elementary proof is left to the reader.

\subsection{Nets of rational functions}

Let $R^d$ be the class of real non-constant rational functions $g$
of degree at most $d$ whose all critical points are
real. 
Consider the full preimage $\gamma=g^{-1}(\bR)$, where
$\bR=\R\cup\{\infty\}$.
This preimage consists of simple analytic arcs which meet only at
critical points. These arcs define a cell decomposition $C(g)$ of the
Riemann sphere $\P^1$ whose $2$-cells (faces)
are components of $\P^1\backslash
\gamma$, $1$-cells (edges) are components of $\gamma\backslash\{
{\mathrm{critical\; points}}\}$ and $0$-cells (vertices) are
the critical points. 
We choose some vertex $v_0$ and call it the
{\em distinguished vertex} of $C(g)$.
The union $\gamma$ of edges and vertices is
the $1$-skeleton of the cell decomposition. Such a cell decomposition
$C=C(g)$ has the following properties:
\vspace{.1in}

\noindent
(i) the $1$-skeleton of $C$ contains $\bR$.

\noindent
(ii) $C$ is symmetric with respect to $\bR$,

\noindent
(iii) all vertices belong to $\bR$
and an even number of edges
meet at each vertex. 
\vspace{.1in}

If $C=C(g)$ the even number in (iii)
is twice the local degree of $g$ at the critical point.
Another important property of our cell decomposition is that
the closure of every cell is homeomorphic to a closed ball
of the same dimension.
It follows that
\vspace{.1in}

\noindent
(iv) no edge can begin and end
at the same point.
\vspace{.1in}

Two cell decompositions $C_1$ and $C_2$ satisfying (i)--(iv)
and having distinguished vertices $v_0^1$ and $v_0^2$
will be called {\em equivalent} if there exists a homeomorphism
$\phi:\P^1\to\P^1$ commuting with reflection with respect
to $\bR$, preserving orientations of $\bR$ and $\P^1$,
mapping cells of $C_1$ onto cells of $C_2$ and $v_0^1$ to $v_0^2.$
An equivalence class of cell decompositions satisfying (i)--(iv)
will be called a {\em net}. The number of faces of a net is even,
we denote it by $2d$ and call the positive integer $d$
{\em the degree of the net.} If $C=C(g)$ then $\deg C=\deg g$.
A net of degree
$d$ has $2d-2$ edges disjoint from the real axis.
Using the Uniformization theorem, we proved in \cite{Ann,CMFT}
that each net comes from a
rational function, and the critical
points of this rational function can be arbitrarily prescribed,
but we do not use this result here, and in fact it will be
deduced from our theorems 1 and 2 in the end of Section 1.

We need the following elementary
\vspace{.1in}

\noindent
{\bf Proposition 2}. {\em Let $(p_j)$ be a convergent
sequence in $G_\R$, and $p=\lim p_j.$  Let $g_j=r(p_j)$ be the corresponding
sequence of rational functions. Then the sets $\gamma_j=g_j^{-1}(\bR)$
converge in the Hausdorff metric to the set $\gamma=g^{-1}(\bR)$,
where $g=r(p)$.}
\vspace{.1in}

This is a simple corollary of Proposition 1,
and the details are left to the reader. 

\vspace{.1in}

\noindent {\bf Corollary 1.}
{\em Suppose that $g_t: t\in[0,1]$ is a continuous path in $R^d$
and each $g_t$ 
has $2d-2$ simple critical points. Let one of these critical
points be $v_0(t)$, a continuous function of $t$.
Then the net of $g_t$ with distinguished point $v_0(t)$
is independent of $t$.}
\vspace{.1in}


\noindent
{\bf Corollary 2.} {\em Let $p_t$ be a continuous path in the Grassmannian
$G$, pa\-ra\-met\-ri\-zed by $[0,1]$. Suppose that 
$g_t=r(p_t)$ belong to $R^d$ and have critical points
$x_0(t),\ldots, x_n(t)$, such that $x_j(t)\neq x_i(t)$ for $0\leq j<i\leq n
$ and
$0\leq t<1$, while for $t=1$ we have $x_0(1)=x_1(1)$ and 
$x_j(t)\neq x_i(t)$ for $1\leq j<i\leq n$. 

Then the degrees of $g_t$ are equal for $t\in[0,1)$, and the degree of
$g_1$ is less than the degree of $g_t,\; t\in[0,1)$
if and only if the
net of $g_0$ contains an edge from $x_0$ to $x_1$.}
\vspace{.1in}

{\em Proof.} According to Proposition 2, the sets $\gamma_t=
g^{-1}_t(\bR)$ vary continuously in the Hausdorff metric.
If there is an edge connecting $x_0(t)$
and $x_1(t)$, the limit of this edge as $t\to 1$ cannot be a
loop because of the property (iv) of the nets.
So the limit belongs
to $\bR$ and thus the limit cell decomposition has fewer faces
than the cell decomposition $C(g_1)$.

In the opposite direction, if this limit has fewer faces, some edge has
to disappear in the limit, and this can only be an edge from $x_0$
to $x_1$. 
\hfill$\Box$

\subsection{Thorns}

Let $X^{2d-2}\subset\Poly_\R^{2d-2}$
be the subset consisting of polynomials
whose all roots are real. Then $X^{2d-2}$ has non-empty interior.  
Here we construct an  open subset of $X^{2d-2}$
such that for every
polynomial $p$ in this subset, the full preimage
$W^{-1}(p)$ consists of $u_d$
distinct real points in $G$. The existence of such a subset
was established by Sottile \cite{Sot1}, but we give a more precise
description of this set following \cite{DCG}.

We fix $d\geq 2$. Consider pairs of integers 
\begin{equation}
\label{rule}
0\leq k_1<k_2\leq d
\end{equation}
and pairs $\q=(q_1,q_2)$ of real polynomials
\begin{equation}
\label{canon}
\begin{array}{rclrl}
q_1(z)&=&&z^{d-1}&+a_{1,d-2}z^{d-2}+\ldots+a_{1,k_1}z^{k_1},\\
\\
q_2(z)&=&z^{d}+&a_{2,d-1}z^{d-1}&+\ldots+a_{2,k_2}z^{k_2}.
\end{array}
\end{equation}
Suppose that all coefficients $a_{i,j}$ are strictly positive, all roots of the
Wronskian determinant $W(\q)=W(q_1,q_2)$ belong to the semi-open interval
\newline
$(-1,0]\subset\R$, and those roots on the open interval $(-1,0)$
are simple. The set of all such polynomial pairs (\ref{canon})
will be denoted by
$b(k_1,k_2)$. The greatest common factor of $\{ q_1,q_2\}$ is $z^{k_1}.$

It is easy to see that the representation
of a point of $G_\R$ by a pair in $b(k_1,k_2)$ is unique.
Setting $k_1=0$ and $k_2=1$ we obtain an open subset
$b(0,1)$ of $G_\R$.

The multiplicity of the root of $W(\q)$
at $0$ is $k=k_1+k_2-1$. We enumerate the negative roots of $W(\q)$ as
\begin{equation}
\label{enumera}
-x_{2d-2}<-x_{2d-3}<\ldots<-x_{k+1}<0.
\end{equation}

Let $\epsilon$ be a positive increasing function on $[0,1]$
satisfying $\epsilon(0)=0, \; \epsilon(x)<x,$ for $x\in(0,1]$.
The set of all such functions will be denoted by $E$.
A {\em thorn} $T$ of dimension
$n$ is a region in $\R^n$ of the form
$$T(n,\epsilon)=
\{(y_1,\ldots,y_n): 0<y_n<\epsilon(1),\; 0<y_{k}<\epsilon(y_{k+1}),
\; 1\leq k\leq n-1\}.$$
Let $w(k,T)$ be the set of polynomials of the form
\begin{equation}
\label{p}
p(x)=x^k(x+x_{2d-2})(x+x_{2d-3})\ldots(x+x_{k+1}),
\end{equation}
where the vectors $(x_{2d-2},\ldots, x_{k+1})$ belong
to a thorn $T$ of dimension $2d-2-k$.


The set $b(d-1,d)$ consists of a single pair $q_1=z^{d-1},\; q_2=z^d$.

Consider the following two operations 
$F^i,\; i=1,2$.
For each pair $(q_1,q_2)$, of the form (\ref{canon}) operation $F^i$
adds to the polynomial $q_i$ one term $az^{k_i-1}$, where
$a>0$ is a small parameter, and leaves the other polynomial of the
pair unchanged. So each operation $F^i$ increases the total number
of non-zero coefficients of a polynomial pair (\ref{canon}) by one:
\begin{equation}
\label{operators}
F^i_a:b(\k)\to b(\k-\e_i),\quad F^i_a(\q)=\q+az^{k_i-1}\e_i,
\end{equation}
where $\k=(k_1,k_2)$ and $(\e_1,\e_2)$ is the standard basis in $\R^2$,
and $a>0$ is a small parameter whose range may depend on $\q$.  
The following rule has to be observed:
\vspace{.1in}

\noindent
{\bf Rule}. {\em Operation $F^i$ is permitted on $b(k_1,k_2)$ if and
only if the outcome of this operation does not violate inequalities}
(\ref{rule}).
\vspace{.1in}

In other words, operation $F^1$ is permitted on $b(k_1,k_2)$ if $k_1>0$,
and operation $F^2$ is permitted on $b(k_1,k_2)$ if $k_2>k_1+1.$

The following result, which is a part of \cite[Proposition 8]{DCG}
shows, among other things, that the $F^i$ are well defined, that is their
result indeed belongs to some
$b(k_1^*,k_2^*)$ for sufficiently small values
of the parameter $a$.
\vspace{.1in}

\noindent
{\bf Proposition 3.} {\em Suppose that
$\k=(k_1,k_2)$ and $i\in\{1,2\}$ satisfy the Rule above.
Suppose that for some thorn $T$ of dimension $2d-2-k$ a set $U\subset b(k_1,k_2)$ is given,
such that
the map $\q\mapsto W(\q): U\to w(k,T)$ is surjective.
Then there exists a thorn $T^*$ of dimension $2d-1-k$
and a set $U^*\subset b(\k-\e_i)$,
such that every $\q^*\in U^*$ has the form $F^i_a(\q),\; \q\in U$, where
$F^i_a$ is defined in $(\ref{operators})$, and $a>0$, 
and the map
\begin{equation}
\label{map}
\q^*\mapsto W(\q^*): U^*\to w(k-1,T^*)
\end{equation}
is surjective.}
\vspace{.1in}

{\em Proof.} We follow \cite[Section 2]{DCG}. First we state three 
elementary lemmas about thorns.
\vspace{.1in}

\noindent
{\bf Lemma 1.} {\em Intersection of any finite set of thorns
of same dimension contains a thorn of the same dimension.}
\vspace{.1in}

{\em Proof.} Take the minimum of their defining functions.
\hfill$\Box$
\vspace{.1in}

\noindent
{\bf Lemma 2.} {\em
Let $T=T(n,\epsilon)$ be a thorn of dimension $n$ in $\R^n=
\{(x_{1},\ldots,x_{n})\}$,
and $U$ its neighborhood
in $\R^{n+1}=\{(x_0,\ldots,x_{n})\}.$
Then $U^+=U\cap\R^{n+1}_{>0}$ contains a thorn $T(n+1,\epsilon_1)$.
}
\vspace{.1in}

{\em Proof.}
There exists a continuous function $\delta_0:T\to\R_{>0}$,
such that $U^+$ contains the set $\{(x_0,\x):\x\in T, 0<x_0<\delta_0(\x)\}$.
Let $\delta(t)$ be the minimum of
$\delta_0$ on the compact subset
$\{\x\in \overline{T(n,\epsilon/2)}:x_{1}
\geq t
\}$
of $T$. Then there exists $\epsilon_0\in E$ with the property $\epsilon_0<
\delta$.
If we define
$\epsilon_1=\min\{\epsilon/2,\epsilon_0\}$, then $T(n+1,\epsilon_1)
\subset U^+$.
\hfill$\Box$
\vspace{.1in}

\noindent
{\bf Lemma 3.} {\em
Let $T=T(n+1,\epsilon)$ be a thorn of dimension $n+1$,
and\newline
$h:T\to\R^{n+1}_{>0},$
$(x_0,\x)\mapsto(y_0(x_0,\x),\y(x_0,\x))$,
a continuous map with the properties:
for every $\x$ such that $(x_0,\x)\in T$
for some $x_0>0$,
the function $x_0\mapsto y_0(x_0,\x)$ is increasing,
and $\lim_{x_0\to 0}\y(x_0,\x)=\x$.
Then the image $h(T)$ contains a thorn.
}

\vspace{.1in}

{\em Proof.} We consider the region $D\in\R^{n+1}$ consisting of $T$,
its reflection $T'$ in the hyperplane $x_0=0$ and the interior
with respect to
this hyperplane
of the common boundary of $T$ and $T'$.
The map $h$ extends to $T'$ by symmetry:
$h(-x_0,\x)=-h(x_0,\x), (x_0,\x)\in T$, and then to the whole $D$
by continuity. It is easy to see that the image of the extended map
contains a neighborhood $U$ of the
intersection of $D$ with the hyperplane $x_0=0$.
This intersection is a thorn $T_1$ in $\R^n=\{(x_0,\x)\in\R^{n+1}:
x_0=0\}$.
Applying Lemma 2
to this thorn $T_1$, we conclude that $U^+$ contains a thorn.
\hfill$\Box$
\vspace{.1in}

We continue the proof of Proposition 3.

Let us fix $\q\in U$, and put $W=W_\q$.
As $W\in w(k,\epsilon),$ we have $\ord W=k$,
where $\ord$ denotes the multiplicity of a root at $0$.
Let $cz^k$ be the term of the smallest degree in $W(z)$.
Then $c>0$, because all roots of $W$ are non-positive.
In fact,
\begin{equation}
\label{c}
c=(k_2-k_1)a_{2,k_2}a_{1,k_1}>0.
\end{equation}
We fix $i\in\{1,2\}$ satisfying the Rule above and define $W^*=W_{\q^*},$ where $\q^*=F_a^i(\q)$.
Then $\ord W^*=k-1$ and the term of the smallest degree
in $W^*(z)$ is $c^*z^{k-1}$, where
\begin{equation}
\label{cstar}
c^*=a(k^*_2-k^*_1)a_{3-i,k_{3-i}}>0,
\end{equation}
We conclude that when $a$ is small enough (depending on $\q$), the
Wronskian $W^*$ has one simple root in a neighborhood of each negative root
of $W$, and in addition, one simple negative root close to zero, and
a root of multiplicity $k-1$ at $0$.
To make this more precise, we denote the negative roots of $W$ and $W^*$ by
\begin{equation}
\label{igrek1}
-x_{2d-2}<\ldots<-x_{k+1}\quad\mbox{and}\quad
-y_{2d-2}<\ldots<-y_{k+1}<-y_k,
\end{equation}
where $y_j=y_j(a)$.
We have
\begin{equation}
\label{trah}
y_j(0)=x_j,\quad\mbox{for}\quad 1\leq j\leq n,\quad\mbox{and}\quad y_k(0)=0.
\end{equation}
Furthermore, if $a$ is small enough (depending on $\q$)
\begin{equation}
\label{trah2}
a\mapsto y_k(a)\quad\mbox{is increasing and continuous}.
\end{equation}
The set $w(k,\epsilon)$ is parametrized by a thorn
$T=T(2d-2-k,\epsilon)$,
where $\x=(x_{k+1},\ldots,x_{2d-2})$. 
There exists a continuous function $\delta:T\to\R_{>0}$, such that
\begin{equation}
\label{11}
\q^*\in b(\k^*),\quad\mbox{for}\quad a\in (0,\delta(\x)),\quad
\x\in T.
\end{equation}
It remains to achieve (\ref{map}) by modifying the thorn $T$.
Consider the set
\begin{equation}
\label{newset}
U^*=\{\q^*=F_a(\q_\x):\x\in T,\, a\in(0,\delta(\x))\}\subset b(\k^*),
\end{equation}
where $\q_\x\in U$ is some preimage under $W$ of the
polynomial (\ref{p}) with
\newline $(x_{k+1},\ldots,x_{2d-2})=\x$. Such preimage exists
by assumption of Proposition 3 that the map $\q\mapsto W_\q,\;
U\to w(k,T)$ is surjective.
We
apply Lemma 2 to the half-neighborhood (\ref{newset})
of
$T$, with $x_k=a$, to obtain a thorn $T_1(2d-k-1,\epsilon_1)$.
Then we apply Lemma~3 to the map $h:T_1\to\R^{2d-k-1}_{>0}$,
defined
by $y_j=y_j(x_0,\x)$, where $y_j$ are as in (\ref{igrek1}), and $x_k=a$.

This map $h$ satisfies all conditions of Lemma 3 in view of
(\ref{trah}) and (\ref{trah2}). This proves (\ref{map}) and Proposition 3.
\hfill$\Box$
\vspace{.1in}

We begin with the single element of $b(d-1,d)$
and apply operations $F^i$ in some sequence, obeying the Rule above,
while possible. As every step decreases $k$ by $1$, the total number
of steps will be $2d-2$. We describe the sequence of steps
by a sequence $\sigma$ of $1$'s and $2$'s of length $2d-2$.
The number $i$ on the $n$-th place in this
sequence indicates that operator $F^i$ was applied on
$n$-th step.
The Rule above translates to the following characterization
of all possible sequences $\sigma$:

a) The numbers of $1$'s and $2$'s in $\sigma$ are equal.

b) In each initial segment of $\sigma$ the number of $1$'s is
not less than the number of $2$'s.

Such sequences are called {\em ballot sequences} (for two candidates),
see, for example, \cite{Stanley}. The number of ballot sequences
of length $2d-2$ is the Catalan number $u_d$.

Applying $2d-2$ times Proposition 3 according to each ballot
sequence $\sigma$ we obtain in the end
an open subset $U_\sigma\subset b(0,1)$
which is mapped surjectively by the Wronski map onto
$w(0,T)$ for some thorn $T$. As the intersection of any finite set of
thorns of the same
dimension contains a thorn of the same dimension by Lemma 1, we may assume
that this set $w(0,T)$ is the same for all sequences~$\sigma$.

When applying Proposition 3 we can make the range of parameter $a$
as small as desired; using this we can assure that the sequence
of coefficients of the pair (\ref{canon}) is monotone: the coefficients
decrease in the order of their appearance. This implies that the
sets $U_\sigma$ with different $\sigma$ are disjoint.

So we obtain $u_d$ disjoint open sets $U_\sigma$, and each of them
is mapped onto $w(0,T)$ continuously and surjectively by the Wronski map $W$.
As the number of preimages of any point under $W$
is at most $u_d$ we conclude that the maps $W:U_\sigma\to w(0,T)$
are homeomorphisms for all $\sigma.$

Thus each point of the open subset $w(0,T)\subset X^{2d-2}$ has
exactly $u_d$ preimages under $W$ and all these preimages are real.
Each of these preimages corresponds to an analytic branch of the inverse
$W^{-1}$ on $w(0,T)$. The branches are enumerated by ballot sequences.

\subsection{Completion of the proof}

Let us fix a thorn $T$ such that each 
polynomial in $w(0,T)$ has $u_d$ different real preimages under
the Wronski map, as in the end of the previous section.

To each of these preimages $\q=(q_1,q_2)$ corresponds a
rational function $r(\q)=q_1/q_2$ in $R^d$ with $2d-2$
distinct real critical points which has a net $\gamma(\q)$.
We take the rightmost critical point of these functions as
distinguished vertices of the nets.

We claim that all these $u_d$ nets are different.
To prove the claim we just
show how to determine the net from the ballot sequence and vice versa.
\vspace{.1in}

\noindent
{\bf Proposition 4.} {\em Let $k=0$, and let $p$ be 
a polynomial in $w(0,T)$ of the form $(\ref{p})$. 
Let $\q=(q_1,q_2)$ be a polynomial pair
as in $(\ref{canon})$ corresponding to a point in $W^{-1}(p)$,
with the ballot sequence $\sigma$,
and $g=r(\q)=q_1/q_2.$

Then the net of $g$ contains an edge 
between $x_m$ and $x_{m+1}$ if and only if the $m$-th member
of the sequence $\sigma$ is $1$.}
\vspace{.1in}

{\em Proof.} It is enough to investigate what happens to a net
when an operator $F^i$ of Proposition 3 is applied.
We see from (\ref{canon}) and (\ref{operators}) that the degree
of $q_1/q_2$ increases if and only if $i=1$. Corollary 2 of
Proposition 2 says that this happens if and only if the net has
an edge between $x_m$ and $x_{m+1}$.
This proves Proposition 4.
\hfill$\Box$
\vspace{.1in}

So we obtained a polynomial $p_0\in w(0,T)$ of degree $2d-2$ with
$2d-2$ real roots whose preimage under
the Wronski map consists of $u_d$ pairs with different nets.
Each of these preimages corresponds to a holomorphic
inverse branch of the Wronski map
in $w(0,T)$. Let $p_1$ be any real
polynomial of degree $2d-2$ with $2d-2$
real roots. Then there exists a path $p_t:t\in[0,1]$ in $X^{2d-2}$
connecting $p_0$ and $p_1$, such that all $p_t$ are
polynomials with $2d-2$ distinct roots. For example one
can connect the corresponding roots of $p_0$ and $p_1$ linearly.
We do analytic continuation of
all the inverse branches of the Wronski map along this path.
As critical points of our rational functions cannot
collide (because the zeros of their Wronskians $p_t$ do not
collide), their nets do not change under the continuation.
Suppose that this analytic continuation to $t=1$ is impossible.
Let $t_0$ be the smallest singular point. Then $p_{t_0}$
is a ramified value of the Wronski map and the full
preimage $W^{-1}(p_{t_0})$ consists of fewer than $u_d$ points.
This full preimage still consists of real rational functions
with all critical points real and distinct, 
so the nets are defined
for all elements of this preimage.
This means that at least two one-parametric families of
rational functions with different nets
tend to the same function with $2d-2$ distinct critical points,
which is impossible by Corollary 1 of Proposition 2.

This proves theorems 1 and 2.
\vspace{.1in}

This proof clearly implies that for any given net there exists
a unique class of real rational functions with all critical points real,
and the critical points of these functions can be chosen arbitrarily,
the result which was established in \cite{Ann} with the help
of the Uniformization theorem and rather complicated topological
considerations.
\vspace{.1in}

To prove Theorem 3, we notice that $1$-skeleton of every net
of degree $d$ can be obtained
as the limit of $1$-skeletons of nets with $2d-2$ vertices.
So every net of degree $d$ actually occurs as a net of a real rational
function of the class of degree $d$ with all critical points real.
Counting the nets of degree $d$
with $q$ vertices of degrees $2a_1\ldots,2a_q$ gives the Kostka number
$K_\a$ (see, for example, \cite[Lemma 3]{egsv}. So there are at least
$K_\a$ classes of rational functions of degree $d$ with prescribed
real critical points. On the other hand, Schubert calculus \cite{Sch}
shows that there are at most $K_\a$ classes of rational functions
with any prescribed critical points of multiplicities $a_1,\ldots,a_q$.
This proves Theorem 3.
\vspace{.1in}

\noindent{\bf Corollary.} {\em To each net of degree $d$ corresponds exactly one
class of real rational functions of degree $d$ with prescribed real
critical points.}

\section{Fuchsian equations}

\subsection{Equations with all polynomial solutions}

Suppose that distinct points
$a_1,\ldots,a_n$ in the complex plane are given.
We want to describe the set of equivalence classes of polynomial pairs
$(y_1,y_2)$,
such that 
$$W(y_1,y_2)\sim (z-a_1)\ldots(z-a_n).$$
This time we {\em do not specify in advance}
the degree of polynomials $y_1,y_2$
but it is
easy to see that it is at most $n+1$. As $W(y_1,y_2)$ has only simple zeros,
the polynomials $y_1$ and $y_2$ are co-prime.

To approach this problem, we introduce a special {\em parametrization
of equivalence classes}
of polynomial pairs whose Wronskian has prescribed zeros.
Recall that for any two
linearly independent functions
$y_1$ and $y_2$ one can write a second order
linear differential equation which has these two functions as solutions:
\begin{equation}
\label{ABC}
\left|\begin{array}{ccc}y&y_1&y_2\\
                          y'&y_1^\prime&y_2^\prime\\
                          y''&y_1^{\prime\prime}&y_2^{\prime\prime}\end{array}
\right|=
Ay''+By'+Cy=0,
\end{equation}
where
$$A=W(y_1,y_2),\quad B=-A'\quad\mbox{and}\quad C=W(y_1^\prime,y_2^\prime).$$

If $y_1$ and $y_2$ are polynomials, then $A,B,C$ are also polynomials,
and we have 
\begin{equation}\label{degs}
\deg B\leq \deg A-1,\quad\deg C\leq \deg A-2.
\end{equation}
These conditions are equivalent to regularity of the singular point
at infinity.
Introducing two rational functions $P=B/A$ and $Q=C/A$, we conclude
that $P(\infty)=Q(\infty)=0$. Furthermore, $P=-A'/A$,  and $A$ has only
simple zeros, so all residues of $P$ in $\C$ are equal to $-1$.
Denoting the residues of $Q$ by $x_j$, we obtain
\begin{equation}
\label{00}
P(z)=-\sum_{j=1}^n\frac{1}{z-a_j}\quad\mbox{and}\quad
Q(z)=\sum_{j=1}^n\frac{x_j}{z-a_j}.
\end{equation}
Evidently, $P,Q$ and $x_j$ depend only on the equivalence class
of the polynomial pair $(y_1,y_2)$. 

Now we write the conditions on $x_j$ which express the fact that all solutions
of the differential equation
\begin{equation}
\label{0}
y''+Py'+Qy=0
\end{equation}
are polynomials.

In a neighborhood of a singular point $a_k$, our equation can be written
in the form:
\begin{equation}
\label{1}
(z-a_k)y''+P_k(z)y'+Q_k(z)y=0,
\end{equation}
with
$$P_k(z)=-1+p_k(z-a_k)+O(z-a_k)^2,\quad\mbox{where}\quad p_k=\sum_{j\neq k}
\frac{1}{a_j-a_k},$$
and 
$$Q_k(z)=x_k+q_k(z-a_k)+O(z-a_k)^2,\quad\mbox{where}\quad
q_k=-\sum_{j\neq k}\frac{x_j}{a_j-a_k}.$$
Our differential equation (\ref{1}) has two linearly independent polynomial
solutions without a common factor, so it has a polynomial
solution of the form
$$y(z)=1+c_1(z-a_k)+c_2(z-a_k)^2+O(z-a_k)^3.$$ 
Differentiating this, we obtain
$$y'(z)=c_1+2c_2(z-a_k)+O(z-a_k)^2,$$
$$y''(z)=2c_2+O(z-a_k).$$
Substituting this to (\ref{1}), we obtain
$$c_1=x_k,$$
and
\begin{equation}
\label{4}
p_kc_1+x_kc_1+q_k=0,
\end{equation}
so
$$x_k^2=-p_kx_k-q_k.$$
Recalling the expressions for $p_k,q_k$, we obtain the following
{\em necessary} condition for the equation (\ref{1}) to have two linearly independent
polynomial solutions:
\begin{equation}
\label{2}
x_k^2=\sum_{j\neq k}\frac{x_j-x_k}{a_j-a_k},\quad k=1,\ldots,n.
\end{equation}

\noindent
{\bf Proposition 5.}
{\em Condition $(\ref{2})$ is necessary and sufficient for a differential
equation $(\ref{0})$ with coefficients 
$(\ref{00})$ to have two linearly independent polynomial solutions.}
\vspace{.1in}

{\em Proof}. It remains to prove sufficiency.
{}From (\ref{1}) we conclude that
all singular points in $\C$ are regular, with exponents
$0$ and $2$. Condition (\ref{4}) guarantees that there is a power
series solution corresponding to the smaller exponent. This implies
that there are two linearly independent holomorphic solutions
in a neighborhood of each singular point. Thus all solutions are
entire functions. By a theorem of Halphen which can be found in
\cite[15.5]{Ince}, if the general solution
of the equation
(\ref{ABC}), where $A,B,C$ satisfy $\deg B\leq\deg A,\;\deg C\leq\deg A$,
is a meromorphic function
in $\C$,
then this meromorphic function has to be of the form
$$\sum R_j(z)e^{\lambda_jz},$$
where $R_j$ are rational functions. But asymptotics at infinity shows that
under the stronger condition (\ref{degs}) on the degrees of the coefficients,
the exponentials cannot be present.
So the general solution is an entire
rational
function, that is a polynomial.
\hfill$\Box$
\vspace{.1in}

\noindent
{\bf Proposition 6.} {\em Every solution of the system $(\ref{2})$
has the following properties:
\begin{equation}
\label{a}
\sum_{k=1}^n x_k=0,
\end{equation}
\begin{equation}
\label{b}
{(n+1)^2-4\sum_{k=1}^nx_ka_k}=s^2,
\end{equation}
where $s$ is an integer such that $n+s$ is odd, and $1\leq s\leq n+1$.
}
\vspace{.1in}

This integer $s$ is the local degree at infinity of the rational function
which is the ratio of two linearly independent solutions of (\ref{0}) with 
$P$ and $Q$ as in (\ref{00}).

{\em Proof of Proposition 6}. By Proposition 5, for every solution $x_1,\ldots,x_n$
of (\ref{2}), all solutions of the differential equation (\ref{0})
are polynomials. Suppose that
$$P(z)=pz^{-1}+O(z^{-2}),\quad Q(z)=qz^{-1}+q^*z^{-2}+O(z^{-3})\quad z\to\infty.$$
Then
$$p=-n,\quad q=\sum_{k=1}^nx_k\quad\mbox{and}\quad q^*=\sum_{k=1}^nx_ka_k.$$
Substituting into (\ref{0}) a polynomial
$$y(z)=z^d+\ldots,$$
we first obtain
$$
qz^{d-1}+O(z^{d-2})=0,\quad z\to\infty,$$
so $q=0$ which proves (\ref{a}). Then we obtain
$$d(d-1)z^{d-2}+dpz^{d-2}+q^*z^{d-2}+O(z^{d-3})=0,$$
so 
$$d^2-d(n+1)+q^*=0.$$
This equation has two solutions, the possible degrees of polynomials:
$$d_{1,2}=\frac{n+1\pm\sqrt{(n+1)^2-4q^*}}{2}.$$
As all solutions of our differential equation are polynomials, there
are two solutions of different degrees.
This implies that both $d_1$ and $d_2$ are
integers
which implies (\ref{b}). Notice that $d_1-d_2=s$.

Rational function $y_1/y_2$ is locally $s$ to $1$ at
infinity, from which the inequality $1\leq s\leq n+1$ follows.
\hfill$\Box$
\vspace{.1in}

So our theorems 1 and 2 are can be restated as
\vspace{.1in}

\noindent
{\bf Theorem 4.} {\em If all $a_k$ are real then all solutions of
$(\ref{2})$ are real. Moreover, each solution $(x_1,\ldots,x_n)$
can be analytically continued as a
function of parameters $(a_1,\ldots,a_n)$ in the region
$a_1<a_2<\ldots<a_n$.}
\vspace{.1in}

Solutions $(x_1,\ldots,x_n)$ correspond to classes of rational functions
of all degrees $d\in [n/2+1,n+1]$ if $n$ is even and
$d\in[(n+1)/2+1,n+1]$ if $n$ is odd, 
having simple critical points at $a_1,\ldots,a_n$
and possibly a critical point at infinity. The order of the critical
point at infinity is $s-1$ where $s$ is defined in (\ref{b}).

System (\ref{2}) has a trivial solution $x_1=\ldots=x_n=0$
which corresponds to the polynomial of degree $n+1$
with critical points $a_1,\ldots,
a_n$. The opposite case is that $s=1$ (then $n$ is even) and
we have rational functions with $n$ prescribed simple critical
points. This case is characterized by the condition
$$q^*=(n^2+2n)/4.$$

System (\ref{2}) can be easily generalized to
$p$-tuples of linearly independent
polynomials satisfying a Fuchsian differential equation
of order $p$. An analog of Theorem 4 for this case would be
equivalent to the Shapiro conjecture for $p$ polynomials. 
For example, for $p=3$ one obtains the following system
of equations with respect to $x_k$ and $u_k$:
\begin{eqnarray*}\displaystyle
x_k^2&=&\sum_{j\neq k}\frac{x_j-x_k}{a_j-a_k}-u_k,\\
x_ku_k&=&\sum_{j\neq k}\frac{u_j-u_k}{a_j-a_k}.
\end{eqnarray*}


\subsection{Equilibria of electric charges in the plane}

Following Stieltjes \cite{Stiltjes}, we state an extremal problem
of potential theory which is equivalent to (\ref{2}).

Suppose as above that 
\begin{equation}
\label{A}
A(z)=(z-a_1)\cdots(z-a_n),\quad\mbox{and}\quad B=-A'
\end{equation}
in (\ref{ABC}), and let us look for polynomials $C$ such that (\ref{ABC})
has a polynomial solution
\begin{equation}
\label{igrek}
y(z)=(z-z_1)\ldots(z-z_m)\quad\mbox{with}\quad z_j\in\C\backslash\{ a_1,\ldots,a_n\}.
\end{equation}
We substitute (\ref{igrek}) in (\ref{ABC}) and put $z=z_k$:
$$A(z_k)y^{\prime\prime}(z_k)-A'(z_k)y'(z_k)=0,$$
or
$$y^{\prime\prime}(z_k)/y'(z_k)-A'(z_k)/A(z_k)=0,$$
which is equivalent to
\begin{equation}
\label{stiltjes}
2\sum_{j\neq k}\frac{1}{z_k-z_j}-\sum_{j=1}^n\frac{1}{z_k-a_j}=0.
\end{equation}
Equation (\ref{stiltjes}) has the following physical interpretation.
Positive unit charges $+1$ are fixed at the points $a_k$,
and negatively charged particles of charge $-2$ each
at the points $z_k$ are allowed
to move in the plane. The particles are repelled or attracted
according to their charges, and the force is
inverse proportional to the distance\footnote{The particles in the plane
can be imagined as infinite uniformly charged wires perpendicular
to the complex plane and interacting by the Coulomb law.}.
So the force between $z$ and $w$
is $c/(\overline{z-w})$, where the bar stands for the complex
conjugation and $c$ is a real constant depending on the charges.
Then the equilibrium condition is expressed by the system
of equations (\ref{stiltjes}).
The energy of such configuration is
$$\displaystyle
E=\log\frac{\prod_{(j,k):j\neq k}|z_k-z_j|^2}{\prod_{(j,k)}|z_k-a_j|}.
$$
The function $U(\z,\a)$ under the logarithm is called the master
function in \cite{Varsch,MV,Sch,Sch2}. Here $\z=(z_1,\ldots,z_m)$
is the variable and $\a=(a_1,\ldots,a_n)$ the parameter.
Equilibrium configurations are the critical points of the master function.
It is evident from physical considerations and easy to prove
that all these equilibria are unstable.

So for a given polynomial $A$ as in (\ref{A}), an equation of the form
\begin{equation}
\label{special}
Ay^{\prime\prime}-A'y'+Cy=0
\end{equation}
has a non-trivial polynomial solution $y$ if and only if 
$$C=(-Ay^{\prime\prime}+A'y')/y,$$
with $y$ given by (\ref{igrek}) and $(z_1,\ldots,z_m)$ satisfying
(\ref{stiltjes}).
\vspace{.1in}

\noindent
{\bf Lemma 4.} {\em If the equation $(\ref{special})$ has one
non-trivial polynomial solution, then all its solutions are
polynomials.}
\vspace{.1in}

This is \cite[Lemma 7]{Varsch}. We include a proof for the reader's
convenience.

{\em Proof.} We apply the usual method of finding a second 
linearly independent solution when one, say $y_1$, is given. We obtain
\begin{equation}\label{expr}
y_2=y_1\int Ay_1^{-2}.
\end{equation}
This can have ramification points only at the zeros of $y_1$.
But as a solution of a differential equation (\ref{special}),
$y_2$ can only have ramification points at the zeros of $A$.
As these two sets are disjoint, $y_2$ is an entire function.
Now, using the expression (\ref{expr}) we conclude that $y_2$
is a polynomial. So all solutions are polynomials.
\hfill$\Box$
\vspace{.1in}

{\em Remark.} An alternative way to derive (\ref{stiltjes}) is
to write the condition that all residues of the integral in (\ref{expr})
are equal to zero.
\vspace{.1in}

We conclude that some solutions of (\ref{stiltjes})
are not isolated; they occur in complex $1$-dimensional families
corresponding to configurations of zeros of all polynomial
solutions of (\ref{special}). However, each two-dimensional space
of polynomials contains a one-dimensional subspace consisting
of polynomials of smaller degree than generic polynomials in
this space.
Such polynomial of the smallest degree gives an isolated
solution of (\ref{stiltjes}). If $n=2d-2$ and $m\leq d-1$, we see
from (\ref{expr})  that $\deg y_2>\deg y_1$, and we conclude
that the critical points of the master function $E$ are all
isolated in this case. These facts about critical points
of the master function and their relation to Fuchsian
equations with polynomial solutions were discovered in \cite{Varsch}.
They can be generalized to Fuchsian equation of arbitrary order.
The equation (\ref{stiltjes}) is a special case of
the ``Bethe Ansatz equation'' in \cite{Varsch}.
 
Suppose now that (\ref{special})
has a polynomial solution. Then all solutions
are polynomials and the ratio of any two linearly independent solutions
is a rational function with all critical points real.
Such a rational function is equivalent to a real rational function
by Theorem 1, so we can find two real linearly independent
solutions. Then $C$ is a real polynomial, and the solution of smallest
degree is proportional to a real polynomial.
Thus Theorem 1 has the following consequence: 
\vspace{.1in}

\noindent
{\bf Theorem 5.} {\em For given real points $a_1,\ldots,a_n$,
each isolated equilibrium
configuration $z_1,\ldots,z_m$ in $(\ref{stiltjes})$
is symmetric with respect to the real line.}
\vspace{.1in}

It is also easy to deduce Theorem 1 from Theorem 5.

\vspace{.1in}

{\em Purdue University, West Lafayette, IN 47907

eremenko@math.purdue.edu

agabriel@math.purdue.edu}

\begin{thebibliography}{1}
\bibitem{Ann} A. Eremenko and A. Gabrielov,
Rational functions with real
critical points and the
B. and M. Shapiro conjecture in real enumerative
geometry, {\em Ann. Math.} 155 (2002) 105-129.
\bibitem{CMFT} A. Eremenko and A. Gabrielov, 
Wronski map and Grassmannians of real codimension 2 subspaces,
{\em Computational Methods and Function Theory} 1 (2001) 1--25.
\bibitem{DCG} A. Eremenko and A. Gabrielov,
Degrees of real Wronski maps, {\em Discrete and Computational Geom.}
28 (2002) 331--347. 
\bibitem{SIAM} A. Eremenko and A. Gabrielov,
Pole placement by static output feedback for generic linear systems,
{\em SIAM J. on Control and Opt.} 41, 1 (2002) 303--312.
\bibitem{egsv} A. Eremenko, A. Gabrielov, M. Shapiro and A. Vainshtein,
Rational functions and real Schubert calculus,
Proc. AMS (electronically published on July 25, 2005)
\bibitem{Goldberg} L. Goldberg,
Catalan numbers and branched coverings by the Riemann sphere. 
{\em Adv. Math.}  85, 2  (1991)  129--144. 
\bibitem{Ince} E. Ince,  {\em Ordinary Differential Equations,}
Longmans, Green and Co.,
London, 1927.
\bibitem{KS} V. Kharlamov and F. Sottile, Maximally inflected real rational
curves,  {\em Mosc. Math. J.}  3, 3  (2003) 947--987, 1199--1200.
\bibitem{MV} E. Mukhin and A. Varchenko,
Critical points of master functions and flag varieties,
{\em Commun. Contemp. Math.} 6, 1 (2004) 111--163.
\bibitem{RosSot} J. Rosenthal and F. Sottile,
Some remarks on real and complex output feedback,
{\em Systems Control Lett.} 33, 2 (1998) 73--80.
\bibitem{Sch} I. Scherbak, Rational functions with prescribed
critical points,  {\em Geom. Funct. Anal.}  12, 6  (2002) 1365--1380.
\bibitem{Sch2} I. Scherbak, Intersections of Schubert varieties and
critical points of the generating function, 
{\em J. London Math. Soc.} (2)  70, 3  (2004) 625--642.
\bibitem{SeS} V. Sedykh and B. Shapiro, On two conjectures concerning convex curves,
to appear in Intl. J. Math.
\bibitem{Sot1} F. Sottile,  The special Schubert calculus is real, 
{\em Electron. Res. Announc. Amer. Math. Soc.} 5 (1999), 35--39. 
\bibitem{Soterag} F. Sottile, Enumerative real algebraic geometry,
in: S. Basu and L. Gonzalez-Vega,
{\em Algorithmic and quantitative real algebraic geometry,}
139--179,
Amer. Math. Soc., Providence, RI, 2003. 
\bibitem{SotShap} F. Sottile,
Real Schubert calculus: polynomial systems
and a conjecture of Shapiro and Shapiro,  
{\em Experiment. Math.} 9, 2 (2000) 161--182.
\bibitem{Stanley} R. Stanley,
{\em Enumerative combinatorics. Vol. 2,}
Cambridge University Press, Cambridge, 1999.
\bibitem{Stiltjes} T. Stieltjes,
Sur certains polyn\^omes qui v\'erifient
une \'equation differentielle lin\'eaire du second ordre et sur la th\'eorie
des fonctions de Lam\'e, {\em Acta Math.} 6 (1885) 321--326;
{\em Oeuvres compl\'etes -- Vol. 1}, Springer, Berlin, 1993, 434--439. 
\bibitem{Varsch} A. Varchenko and I. Scherbak,
Critical points of functions,
$\mathfrak s\mathfrak l\sb 2$ representations,
and Fuchsian differential equations
with only univalued solutions.
{\em Mosc. Math. J.}  3, 2  (2003) 621--645, 745.
\end{thebibliography}
\end{document}